\documentclass{birkjour}
\usepackage{amssymb}
\usepackage{amsmath}
\usepackage{mathrsfs}
\usepackage{amsfonts}
\usepackage{latexsym}
\usepackage{color}
\usepackage{soul,cancel}

\newtheorem{theorem}{Theorem}
\newtheorem{obs}[theorem]{Observation}
\newtheorem{prop}[theorem]{Proposition}
\newtheorem{cor}[theorem]{Corollary}
\newtheorem{lem}[theorem]{Lemma}

\usepackage{rawfonts}

\newfont{\OOO}{cmr10 scaled 2986}
\newfont{\OO}{cmr10 scaled 1440}

    \usepackage[OT2,OT1]{fontenc}
    \newcommand\cyr{%
    \renewcommand\rmdefault{wncyr}%
    \renewcommand\sfdefault{wncyss}%
    \renewcommand\encodingdefault{OT2}%
    \normalfont
    \selectfont}
    \DeclareTextFontCommand{\textcyr}{\cyr}

\newcommand{\C}{\mathbb{C}}

\newcommand{\bN}{{\mathbf N}}

\newcommand{\bU}{{\mathbf U}}

\newcommand{\bW}{{\mathbf W}}

\newcommand{\bY}{{\mathbf M}}

\newcommand{\cA}{{\mathcal A}}

\newcommand{\cH}{{\mathcal H}}
\newcommand{\cB}{{\mathcal B}}

\newcommand{\cM}{{\mathcal M}}

\newcommand{\la}{\lambda}

\newcommand{\si}{\sigma}

\newcommand{\n}{\|}

\newcommand{\Ke}{{\rm Ker}\,}
\newcommand{\Ra}{{\rm Ran}\,}

\newcommand{\eq} [1] {\begin{equation}\label{#1}\quad}
\newcommand{\en} {\end{equation}}

\newcommand{\abs}[1]{\left\vert#1\right\vert}

\newcommand{\diag}{\operatorname{diag}}

\newcommand{\sgn}{\operatorname{sgn}}

\begin{document}
\title[Similarity between two projections]
{Similarity between two projections}

\author[A. B\"ottcher]{Albrecht B\"ottcher}
\address{Fakult\"at f\"ur Mathematik\\
TU Chemnitz\\
09107 Chemnitz\\
Germany}
\email{aboettch@mathematik.tu-chemnitz.de}

\author[B. Simon,]{Barry Simon,}
\address{Departments of Mathematics and Physics\\
Mathematics 253-37\\
California Institute of Technology\\
Pasadena, CA 91125\\
USA}
\email{bsimon@caltech.edu}
\thanks{Research of the second author supported in part by NSF grants DMS-1265592 and DMS-1665526 and in part by Israeli BSF Grant No. 2014337.
 The third author was supported in part by Faculty Research funding from the Division of Science and Mathematics, New York University Abu Dhabi.}
\author[I. Spitkovsky]{Ilya Spitkovsky}
\address{Division of Science\\
New York  University Abu Dhabi (NYUAD)\\
Saadiyat Island\\
P.O. Box 129188\\
Abu Dhabi\\
UAE}
\email{ims2@nyu.edu, imspitkovsky@gmail.com}

\begin{abstract}
Given two orthogonal projections $P$ and $Q$, we are interested in all unitary operators $U$ such
that $UP=QU$ and $UQ=PU$. Such unitaries $U$ have previously
been constructed by Wang, Du, and Dou and also by one of the authors.
One purpose of this note is to compare these constructions.
Very recently, Dou, Shi, Cui, and Du described all unitaries $U$ with
the required property. Their proof is via the two projections theorem by Halmos.
We here give a proof based on the
supersymmetric approach by Avron, Seiler, and one of the authors.
\end{abstract}


\subjclass{Primary 47A62, Secondary 46H15, 46L89, 47A67, 47C15}

\keywords{Intertwining operators, intertwining unitaries, similar projections, two projections}

\maketitle

\section{Introduction}\label{SI}

Let $\cA$ be an algebra with unit $I$ and let $P,Q \in \cA$ be two idempotents, that is, elements
satisfying $P^2=P$ and $Q^2=Q$. We are interested in invertible elements $V \in \cA$ such that $VP=QV$, or
equivalently,
\begin{equation}
Q=VPV^{-1}. \label{v}
\end{equation}
A stronger question is to find an invertible $V \in \cA$ such that $VP=QV$ and $VQ=PV$, which may also
be written as
\begin{equation}
Q=VPV^{-1}, \quad P=VQV^{-1}. \label{vv}
\end{equation}
The strongest version of the problem is to look for a $V$ satisfying~(\ref{v}) and the equality $V^2=I$,
in which case~(\ref{vv}) is automatically valid. If $\cA$ is an algebra with an involution and $P,Q$ are
selfadjoint, which means that $P=P^*$ and $Q=Q^*$, it is natural to ask for unitary elements $V$, i.e., elements
satisfying $V^{-1}=V^*$, which ensure~(\ref{v}) or~(\ref{vv}).
Finally, in case elements $V$ with the required properties exist, we want to describe all of them.
As we were polishing our paper, we learned that~\cite{Dou} had addressed the same problem with some
but certainly not complete overlap as we will explain below.

We begin with a simple observation.

\begin{obs} \label{Obs 1}
Let $V_0$ be an invertible element of $\cA$ such that $Q=V_0 P V_0^{-1}$ and $P=V_0 Q V_0^{-1}$. Then
an invertible element $V \in \cA$ satisfies the equalities~{\rm (\ref{vv})} if and only if $V=CV_0$ with an invertible
element $C \in \cA$ that commutes with both $P$ and $Q$.
\end{obs}

\medskip
{\em Proof.} If $C$ commutes with $P,Q$, then $CV_0 P=C QV_0=QCV_0$ and $CV_0 Q=CPV_0=PCV_0$, which proves the ``if'' part.
Conversely, if~(\ref{vv}) holds for $V=CV_0$, then $CV_0 P=QCV_0$ and since also $CV_0 P=C QV_0$, we conclude
that $QCV_0=C QV_0$, which implies that $QC=CQ$. Analogously one obtains that $PC=CP$. $\;\:\square$

\medskip
Given two idempotents $P$ and $Q$, we put, following~\cite{Av}, $A=P-Q$ and $B=I-P-Q$.
Obviously, an element in $\cA$ commutes with both $P$ and $Q$ if and only if it commutes with both
$A$ and $B$. We have
\begin{eqnarray*}
& & BP=(I-P-Q)P=-QP=Q(I-P-Q)=QB, \\
& &  BQ=(I-P-Q)Q=-PQ=P(I-P-Q)=PB,
\end{eqnarray*}
that is, the two intertwining relations are in force, and if $B$ would be invertible, we would get~(\ref{vv})
with $V=B$. As observed in~\cite{Av}, we also have
\begin{equation}
A^2+B^2=I, \quad AB+BA=0,\label{Pyth}
\end{equation}
and hence $B$ is invertible if and only if $I-A^2$ is invertible. If $\cA$ is a Banach algebra, the invertibility of $I-A^2$
is guaranteed by the inequality $\n A\n <1$. Thus, we can give some quick answers in this case.

\begin{prop} \label{Prop 1}
Let $\cA$ be a Banach algebra and $P,Q$ be idempotents such that $\n P-Q\n <1$. Then  $V=I-P-Q$  is invertible
and satisfies~{\rm (\ref{vv})}. Moreover, there exists an invertible square root
$(I-(P-Q)^2)^{1/2}$ in $\cA$, and
\[V=(I-(P-Q)^2)^{-1/2}(I-P-Q)\]
satisfies~{\rm (\ref{vv})} along with the equality $V^2=I$. If $\cA$ is even a $C^*$-algebra and $P,Q$ are selfadjoint,
then the $V$ just constructed is selfadjoint and unitary.
\end{prop}

{\em Proof.} Since $\n A\n <1$, the power series for $(1-\la)^{-1/2}$ ($|\la| <1$) gives us an invertible
element $C=(I-A^2)^{-1/2}$ by a power series in $A^2$. The operator $C$ commutes with $A$. Since $A^2$ commutes with $B$,
so also does $C$. We know that the invertibility of $I-A^2$ implies the invertibility of $B$,
so that~(\ref{vv}) is satisfied with $V=B$. Observation~\ref{Obs 1} now shows that $V=CB$ also satisfies~(\ref{vv}).
In addition, $V^2=C^{2}B^2=(I-A^2)^{-1}B^2=I$. Finally, if $P$ and $Q$ are selfadjoint, then $A, B, C, V$ are also
selfadjoint and $V$ is unitary. $\;\:\square$

\medskip
{This proof is essentially Kato's~\cite{Kato1, Kato2}. However, Kato used
\[V=PQ+(I-P)(I-Q), \quad \widetilde{V}=QP+(I-Q)(I-P)\] along with the identity
$V\widetilde{V}=\widetilde{V}V=I-A^2$ to get $VP=QV$ and $\widetilde{V}Q=P\widetilde{V}$,
which gives~(\ref{v}). After multiplying these $V$ and $\widetilde{V}$ by $(I-A^2)^{-1/2}$,
he obtained new $V$ and $\widetilde{V}$, which satisfy the same intertwining relations
and also the equalities $V\widetilde{V}=\widetilde{V}V=I$. Our choice $V=I-P-Q$ results in~(\ref{vv}),
not just~(\ref{v}).}

This note is about what can be said if $\n A\n =\n P-Q\n \ge 1$. At this point we remark that the requirement
$\n A\n <1$ may be replaced by first requiring that $1$ is not in the spectrum $\si(A^2)$ of $A^2$
and by secondly requiring that there is a cut of the complex plane from the origin to infinity that does not meet $1-\sigma(A^2)$.
In that case there is an analytic branch of the function $\lambda \mapsto (1-\lambda)^{1/2}$ in an open neighborhood
of $\sigma(A^2)$ and hence the usual formula
\[(I-A^2)^{1/2}=\frac{1}{2\pi i}\int_\Gamma (1-\lambda)^{1/2} (\lambda-A^2)^{-1}d\lambda\]
gives an invertible operator $(I-A^2)^{1/2}$ that commutes with $A$. The rest of the arguments of the above proof then shows
that $V=(I-A^2)^{-1/2}B$ satisfies~(\ref{vv}) and that,
in addition, $V^2=(I-A^2)^{-1}B^2=I$. Notice that the cut mentioned exists in particular if $A=P-Q$ is compact.
In that case $B$ is Fredholm of index zero and hence invertibility is equivalent to injectivity.

The note is organized as follows. After a short section concerning
two skew projections $P$ and $Q$, we consider the case of two orthogonal projections. In Section~\ref{S Orth}, we cite
existing results and constructions and compare them. A new proof for the description of all unitaries satisfying~(\ref{vv})
is given in Section~\ref{S Super}, and in Section~\ref{S Halmos} we quote the result of~\cite{Dou}.
In the final Section~\ref{S Add} we embark on the questions whether the intertwining unitaries belong
to the $W^*$- and $C^*$-algebras generated by $P$ and $Q$.

\section{Skew projections} \label{SSkew}

Let $\cA$ be a complex Banach algebra with unit $I$ and let $P,Q$ be two idempotents in $\cA$.
We denote by ${\rm alg}(P,Q)$ the smallest closed subalgebra of $\cA$ that contains $I,P,Q$.
The following result reveals that the search for an element $V$ satisfying~(\ref{v}) must
go beyond ${\rm alg}(P,Q)$ if $B=I-P-Q$ is not invertible.

\begin{prop} \label{Prop Skew 1}
The following are equivalent:

\smallskip
{\rm (i)} there exists an element $V \in {\rm alg}(P,Q)$ that is invertible
in ${\mathcal A}$ and satisfies $VPV^{-1}=Q$,

\smallskip
{\rm (ii)} $B=I-P-Q$ is invertible in $\cA$,

\smallskip
{\rm (iii)} $1$ is not in the spectrum of $A^2=(P-Q)^2$,

\smallskip
{\rm (iv)} $P+2Q-I$ and $P+2Q-2I$ are invertible in $\cA$.

\smallskip
\noindent
If {\rm (ii)} holds, then $BPB^{-1}=Q$ and $BQB^{-1}=P$.
\end{prop}

{\em Proof.} The equivalence (ii) $\Leftrightarrow$ (iii) follows from the equality
$A^2+B^2=I$. We also know that (ii) implies $BPB^{-1}=Q$ and $BQB^{-1}=P$ and thus~(i).
All other implications can be shown by employing the two projections theorem of
Gohberg and Krupnik~\cite{GK1,GK2} and Roch and Silbermann~\cite{RS}; see
also~\cite{N,BSpit10, Roch}. Example~6.4 of~\cite{BSpit10} contains the equivalence
(ii) $\Leftrightarrow$ (iv). We are therefore left with the implication (i) $\Rightarrow$ (iv).
So suppose (i) holds. If $P+2Q-I$ is not invertible in $\cA$, then the theorem by Gohberg, Krupnik, Roch,
and Silbermann tells us that there exists an algebra homomorphism $G_1: {\rm alg}(P,Q) \to \C$
such that $G_1(V) \neq 0$, $G_1(P)=1$, and $G_1(Q)=0$. Thus, if $VP=QV$, then $G_1(V)\cdot 1= 0\cdot G_1(V)$,
which is impossible. The same theorem shows that if $P+2Q-2I$ is not invertible in $\cA$,
then there is an algebra homomorphism of ${\rm alg}(P,Q)$ into $\C$ satisfying
$G_2(V)\neq 0$, $G_2(P)=0$, and $G_2(Q)=1$, which gives $G_2(V)\cdot 0 =1 \cdot G_2(V)$,
being again a contradiction. This proves the implication (i) $\Rightarrow$ (iv). $\;\:\square$

\medskip
We now turn to the Hilbert space case. Let $\cH$ be a complex Hilbert space and denote by $\cB(\cH)$
the $C^*$-algebra of all bounded linear operators on the space $\cH$. Let $P, Q \in \cB(\cH)$ be two projections,
$P^2=P$ and $Q^2=Q$.
Clearly, if there is an invertible $V \in \cB(\cH)$ such that~(\ref{v}) holds, then the ranges
$\Ra P$ and $\Ra Q$ as well as the kernels $\Ke P$ and $\Ke Q$ must have the same dimension. This is even sufficient for the existence of
an invertible operator $V$ satisfying $VP=QV$. Indeed, in that case we may choose orthonormal bases $\{e_j\}_{j \in J}$
and $\{e_k\}_{k \in K}$ in $\Ra Q$ and $\Ke Q$ and  $\{f_j\}_{j \in J}$
and $\{f_k\}_{k \in K}$ in $\Ra P$ and $\Ke P$, and the operator $V \in \cB(\cH)$ which
maps $e_i$ to $f_i$ is invertible and satisfies $VP=QV$.

{In general, we cannot achieve (\ref{v}) with a unitary operator $V$ and it is also not
possible to ensure~(\ref{vv}) with an invertible operator $V$. To see this,
consider the two-dimensional case and take
\[P=\left(\begin{array}{rr} 0 & -1 \\ 0 & 1 \end{array}\right), \quad
Q=\left(\begin{array}{rr} 1 & 2 \\ 0 & 0 \end{array}\right).\]
These two idempotent matrices have both the rank one and, in addition, satisfy
\[\dim (\Ra P \cap \Ke Q)=\dim (\Ra Q \cap \Ke P)=\{0\},\]
a condition that will make its debut later.
It is easily seen that $2 \times 2$ matrices $V$ and $\widetilde{V}$ obey $VP=QV$
and $\widetilde{V}Q=P\widetilde{V}$ if and
only if
\[V=\left(\begin{array}{cc} -2c & b \\ \,\,c & c \end{array}\right),\quad
\widetilde{V}=\left(\begin{array}{cc} -c & b \\ c & 2c \end{array}\right).\]
But such matrices are never unitary, and the equality $V = \widetilde{V}$ holds if
and only if $c=0$, and in that case the matrices are not invertible.}

For complementary projections, in which case $B=I-P-Q$ is the zero
operator and thus as far away from invertibility as possible, the equality $VP=QV$ implies
$VQ=PV$. Thus, for complementary projections with equal range dimensions we may guarantee~(\ref{vv}) with an invertible $V$.

\section{Orthogonal projections} \label{S Orth}

In the case where $P,Q \in \cB(\cH)$  are orthogonal projections, $P^2=P=P^*$ and $Q^2=Q=Q^*$, much is known.
It was already in 1947, when Sz.-Nagy addressed the question of finding a unitary $U$ such $UP=QU$.
We mentioned above that in the early 1950s, Kato found such a $U$ under the condition that $\| P-Q\| <1$.
An alternative
derivation of the formula for $U$ suggested by Sz.-Nagy and Kato was given by
Mityagin~\cite{Mit}.
In the paper~\cite{Av}, Avron, Seiler, and one of the authors
observed that if the inequality $\n P-Q\n <1$ holds, then one can even find a unitary $U$ such that $UP=QU$ and $UQ=PU$, that is,
such that
\begin{equation}
Q=UPU^*, \quad P=UQU^*, \label{uu}
\end{equation}
and it was shown that $U=\sgn (I-P-Q)$ is as desired and that
Kato's unitary is just equal to $U(I-2Q)=(I-2P)U$ with this $U$.
We refer to Section~5 of~\cite{SimonKato} for more on the early history of the topic
and for precise references.

Wang, Du, and Dou~\cite{WangDuDou} finally proved that a unitary $U$
satisfying~(\ref{uu}) exists if and only if
\begin{equation}
\dim (\Ra P \cap \Ke Q)=\dim (\Ra Q \cap \Ke P). \label{dimdim}
\end{equation}
The ``if'' part of their proof uses Halmos' two projections theorem, the ``only if'' part is through expressing things in a basis.
A proof of this result based on the identities~(\ref{Pyth}) is in~\cite{Simon17}. Since $AB + BA = 0$ is the signature of supersymmetry,
we will call this the supersymmetric approach. Very recently, Dou, Shi, Cui, and Du~\cite{Dou}
gave a description of all unitaries $U$ with the property~(\ref{uu}). Their proof is entirely via Halmos' two projections theorem.
The purpose of what follows is twofold: we will first uncover the connection between~\cite{Av,Simon17} and~\cite{WangDuDou},
and we will secondly state the result of~\cite{Dou} in slightly modified form and derive it
on the basis of the supersymmetric
approach evoked by~(\ref{Pyth}).

We finally remark that the problem of replacing the original equation $Q=UPU^*$ with the stronger
requirement~(\ref{uu}) is not a purely academic matter. This problem actually came from
trying to understand some relations between pairs of projections that appeared in the analysis
of the quantum Hall effect~\cite{AMS, Hall}.

The two projections theorem of Halmos~\cite{Hal69} (see also \cite{BSpit10,Roch} and references therein) is as follows.
Let $P$ and $Q$ be two orthogonal projections acting on the Hilbert space $\mathcal H$, let
$\mathcal L$ and $\mathcal N$ denote the ranges of $P$ and $Q$, respectively, and put
\eq{m} {\mathcal M}_{01}=
{\mathcal L}\cap{\mathcal N}^\perp, \quad  {\mathcal M}_{10}=
{\mathcal N}\cap{\mathcal L}^\perp. \en
Clearly, ${\mathcal M}_{01}=\Ra P \cap \Ke Q$ and ${\mathcal M}_{10}=\Ra Q \cap \Ke P$,
and hence~(\ref{dimdim}) is the same as the equality
\begin{equation}
\dim \cM_{01} = \dim \cM_{10}. \label{mm}
\end{equation}
In addition to (\ref{m}), let
\[{\mathcal M}_{00}={\mathcal L}\cap{\mathcal N}, \quad {\mathcal M}_{11}={\mathcal L}^\perp\cap{\mathcal N}^\perp\]
and
\[{\mathcal M}={\mathcal L}\ominus({\mathcal M}_{00}\oplus{\mathcal M}_{01}), \quad {\mathcal M}'={\mathcal L}^\perp\ominus({\mathcal M}_{10}\oplus{\mathcal M}_{11}).\]
One can show that $\cM$ and $\cM'$ have the same dimension. Let $W: \cM \to \cM'$ be any unitary operator and put $\bW={\rm diag}[I_\cM, W]$.
Let finally $H$ be the restriction of $PQP$ to $\cM$. By construction, $H$ is a Hermitian operator with its spectrum in $[0,1]$ and $0,1$ not being its eigenvalues.
The operators $P$ and $Q$ can now be represented as
\begin{eqnarray}
& & \hspace{-6mm}P=I_{\mathcal M_{00}}  \oplus  I_{\mathcal M_{01}}  \oplus  0_{\mathcal M_{10}}  \oplus  0_{\mathcal M_{11}}
\oplus  {\mathbf W}^*\!\!\begin{bmatrix} I_{\mathcal M} & 0 \\ 0 & 0_{{\mathcal M}'} \end{bmatrix}\!\!{\mathbf W}, \label{Halm1}\\
& & \hspace{-6mm}Q =  I_{\mathcal M_{00}}  \oplus  0_{\mathcal M_{01}} \oplus  I_{\mathcal M_{10}} \oplus 0_{\mathcal M_{11}}
\oplus  {\mathbf W}^*\!\!\begin{bmatrix}  H & \!\!\!\!\sqrt{H(I\!\!-\!\!H)} \,\\ \sqrt{H(I\!\!-\!\!H)} & \!\!\!\!I\!\!-\!\!H \end{bmatrix}\!\!{\mathbf W}.\label{Halm2}
\end{eqnarray}

Armed with this representation we can compare the results of~\cite{Av,Simon17} and~\cite{WangDuDou}.
First of all, we see that if $\n P-Q\n <1$, then we necessarily have $\cM_{01}=\cM_{10}=\{0\}$.
Wang, Du, and Dou~\cite{WangDuDou} assumed~(\ref{mm}) and showed (in a slightly different notation) that if $S:\mathcal M_{10} \to \mathcal M_{01}$  is an arbitrary unitary
operator, then
\begin{equation}
 \bU= I_{\mathcal M_{00}} \oplus \begin{bmatrix} 0 & S \\ S^* & 0\end{bmatrix} \oplus I_{\mathcal M_{11}} \oplus  {\mathbf W}^*\begin{bmatrix}  \sqrt{H} & \sqrt{I-H}\: \\ \sqrt{I-H} & -\sqrt{H} \end{bmatrix}{\mathbf W}, \label{U}
\end{equation}
is unitary and satisfies \eqref{uu}.
Papers~\cite{Av,Simon17} consider the generic case, that is, the case where $\mathcal M_{ij}=\{0\}$ for $i,j=0,1$. The intertwining operator was chosen as $\sgn(I-P-Q)$.
Note that in the absence of the first four direct summands in~\eqref{Halm1} and~\eqref{Halm2} we have
\begin{equation}
 P= {\mathbf W}^*\begin{bmatrix} I & 0 \\ 0 & 0 \end{bmatrix}{\mathbf W}, \quad Q= {\mathbf W}^*\begin{bmatrix}
 H & \sqrt{H(I-H)}\: \\ \sqrt{H(I-H)} & I-H \end{bmatrix}{\mathbf W}, \label{tpg}
\end{equation}
and so
\begin{eqnarray*}
& & B=I-P-Q={\mathbf W}^*\begin{bmatrix}  -H & -\sqrt{H(I-H)}\: \\ -\sqrt{H(I-H)} & H \end{bmatrix}{\mathbf W}, \\
& & B^*B=B^2={\mathbf W}^*\begin{bmatrix}  H & 0 \\ 0 & H \end{bmatrix}{\mathbf W}, \quad
\abs{B}:=\sqrt{B^2}={\mathbf W}^*\begin{bmatrix}  \sqrt{H} & 0 \\ 0 & \sqrt{H} \end{bmatrix}{\mathbf W}.
\end{eqnarray*}
The equality $B=|B|\sgn(B)$ now implies that
\[ \sgn{(I-P-Q)}=\sgn{(B)}= {\mathbf W}^*\begin{bmatrix} -\sqrt{H} & -\sqrt{I-H} \:\\ -\sqrt{I-H} & \sqrt{H} \end{bmatrix}{\mathbf W}. \]
Up to the sign, this is the same as the last direct summand in \eqref{U}. Also notice that in the case of selfadjoint projections
the $V=(I-A^2)^{-1/2}B$ in Proposition~\ref{Prop 1} is nothing but $(B^2)^{-1/2}B=(B^*B)^{-1/2}B=\sgn (B)$.

\section{Orthogonal projections: the supersymmetric approach} \label{S Super}

From the previous section we know that
\begin{equation}
\bU:=I_{\mathcal M_{00}} \oplus \begin{bmatrix} 0 & S \\ S^* & 0\end{bmatrix} \oplus I_{\mathcal M_{11}} \oplus \sgn(B)\label{buu}
\end{equation}
is unitary and satisfies $\bU P=Q\bU$ and $\bU Q=P\bU$. Moreover, condition~(\ref{mm}) ensures that the selfadjoint operators $A$ and $B$ are injective.
Using the decomposition of the Hilbert space $\cH$ into the positive and negative spectral subspaces of $A$ we obtain that
\[A=\bN\begin{bmatrix}a_1 & 0 \\ 0 & -a_2\end{bmatrix}\bN^*\]
with a unitary operator $\bN$ and with positive definite $a_1,a_2$. Let
\[B=\bN\begin{bmatrix}b_1 & b \\ b^* & b_2\end{bmatrix}\bN^*\]
be the respective block matrix representation of $B$.

\begin{lem} \label{Lem a}
The block representations just introduced are actually of the form
\begin{equation}
A=\bY\begin{bmatrix} a & 0 \\ 0 & -a \end{bmatrix}\bY^*, \quad  B=\bY\begin{bmatrix} 0 & \sqrt{I-a^2}\:\\ \sqrt{I-a^2} & 0 \end{bmatrix}\bY^* \label{can}
\end{equation}
where $\bY$ is a unitary operator and $a$ is a positive definite operator which has its spectrum in $[0,1]$ and does not have $1$ as an eigenvalue.
\end{lem}

{\em Proof.}
The anticommutativity condition $AB+BA=0$ implies in particular that $a_1b_1+b_1a_1=0$. Thus, the operator $a_1b_1$ has zero Hermitian part, and so
it is of the form $a_1b_1=iT$ with a selfadjoint operator $T$. It follows that the spectrum of $a_1b_1$ is purely imaginary.
Since $\si(XY)\cup\{0\}=\si(YX)\cup\{0\}$, the spectrum of $a_1^{1/2}b_1a_1^{1/2}$ also is purely imaginary.
On the other hand, the latter operator is selfadjoint, and hence its spectrum is real. Combining the two statements we see that the selfadjoint operator
$a_1^{1/2}b_1a_1^{1/2}$ has zero spectrum and thus itself is zero. From the injectivity of $a_1$ we conclude that $b_1=0$.
It can be shown similarly that $b_2=0$, and so $B$ simplifies to  \[\bN\begin{bmatrix}0 & b \\ b^* & 0\end{bmatrix}\bN^*.\] Since $B$ is injective,
$b$ also is injective and its range is dense. It therefore admits a polar
representation $b=u|b|$ with a unitary factor $u$. With $\bY=\bN \diag[u,I]$, we may
therefore write
\[A=\bY\begin{bmatrix} c_1 & 0 \\ 0 & -c_2 \end{bmatrix}\bY^*, \quad  B=\bY\begin{bmatrix} 0 & |b|\:\\ |b| & 0 \end{bmatrix}\bY^*.\]
In particular, the positive and negative spectral subspaces of $A$ have the same dimension.
Using now the equality $A^2+B^2=I$ we conclude that $c_1^2+|b|^2=I$  and $c_2^2+|b|^2=I$,
which implies that $c_1^2=c_2^2$, and since $c_1$ and $c_2$ are positive definite, it follows that $c_1=c_2 =:a$. Finally,
again because $A^2+B^2=I$, we obtain that $a^2+|b|^2=I$, that is, $|b|=\sqrt{I-a^2}$.
$\;\:\square$

\begin{theorem} \label{Theo Super}
Let $P,Q \in \cB(\cH)$ be orthogonal projections and suppose~{\rm (\ref{mm})} holds. Then all unitary operators
$U$ satisfying~{\rm (\ref{uu})} are given by
\begin{equation}
U=U_0 \oplus \begin{bmatrix} 0 & U_{10} \\ U_{01} & 0\end{bmatrix}\oplus U_1 \oplus{\bY}\begin{bmatrix} 0 & v \\ v & 0 \end{bmatrix}
\bY^*. \label{Ugen1}
\end{equation}
Here $U_j$, $U_{ij}$ are arbitrary unitary operators acting on ${\mathcal M}_{jj}$ and from ${\mathcal M}_{ji}$ onto ${\mathcal M}_{ij}$, respectively,
$v$ is an arbitrary unitary operator commuting with $a$, and $\bY$ is the unitary operator from Lemma~\ref{Lem a}.
\end{theorem}

{\em Proof.}
Let $U$ be a unitary operator satisfying~(\ref{uu}).
We may write $U=Z\mathbf U$ with a unitary operator $Z$ and $\bU$ being the operator~(\ref{buu})
and may represent $Z$ in the form
\[ Z= \left(\bigoplus_{i,j=0,1}Z_{ij}\right)\oplus \bY V \bY^*\]
with $V=\begin{bmatrix} v_{11} & v_{12} \\ v_{21} & v_{22} \end{bmatrix}$.
From Observation~\ref{Obs 1} we infer that $\bY V\bY^*$ must commute with $P,Q$,
hence with $A,B$, and thus also with $\sgn(A),\sgn(B)$. Lemma~\ref{Lem a} shows that
\begin{equation} \label{sgn}
\sgn( A)=\bY\begin{bmatrix} I & 0 \\ 0 & -I \end{bmatrix}\bY^*, \quad  \sgn (B)=\bY\begin{bmatrix} 0 & I \\ I & 0 \end{bmatrix}\bY^*,
\end{equation}
and so every operator commuting with $\sgn (A),\sgn (B)$ has the form
\begin{equation}
\bY\begin{bmatrix} v & 0 \\ 0 & v \end{bmatrix}\bY^*.\label{mvm}
\end{equation}
For $\bY V \bY^*$ to commute with $A,B$ themselves, $v$ has to commute with $a$.
Thus,
\[ Z= \left(\bigoplus_{i,j=0,1}Z_{ij}\right)\oplus{\bY}\begin{bmatrix} v & 0 \\ 0 & v \end{bmatrix}\bY^*.\]
This operator is unitary if and only if so are $Z_{ij}$ and $v$. Taking into account that
\[{\bY}\begin{bmatrix} v & 0 \\ 0 & v \end{bmatrix}\bY^* \sgn(B)=
{\bY}\begin{bmatrix} v & 0 \\ 0 & v \end{bmatrix}\bY^*
\bY\begin{bmatrix} 0 & I \\ I & 0 \end{bmatrix}\bY^*
= {\bY}\begin{bmatrix} 0 & v \\ v & 0 \end{bmatrix}\bY^*,\]
we see that $U=Z\bU$ has the form~\eqref{Ugen1}.
Conversely, as~(\ref{mvm}) commutes with $A,B$ and thus with $P,Q$, Observation~\ref{Obs 1} implies that an operator $U$ as in the theorem
satisfies~\eqref{uu}. $\;\:\square$

\medskip
We remark that if $A$ (equivalently, $a$) has simple spectrum, then $v$ must be a function of $a$; see, for instance, page 308 of~\cite{SimonOT}.
In this case $A$ can be realized as
the operator $M_x$ of multiplication by the independent variable on $L^2({\mathbb R},d\mu)$ with some measure $\mu$; see page~303 of~\cite{SimonOT}.
Then $a$ is nothing but the multiplication operator by $x$ on $L^2({\mathbb R}_+,d\mu)$, and $v$ is the operator of multiplication by an arbitrary unimodular function
on $L^2({\mathbb R}_+,d\mu)$.

\section{Orthogonal projections: the result of Dou, Shi, Cui, and Du} \label{S Halmos}

In slightly modified notation, the result of~\cite{Dou} is as follows.
For convenience, we include a proof.

\begin{theorem}[Dou, Shi, Cui, Du]\label{th:genform}
Let $P,Q$ have the representations~\eqref{Halm1} and~\eqref{Halm2} and satisfy~\eqref{mm}. Then all unitary
operators $U$ satisfying \eqref{uu} are given by the formula
\begin{equation}
  U=U_0 \oplus \begin{bmatrix} 0 & U_{10} \\ U_{01} & 0\end{bmatrix}\oplus U_1 \oplus{\mathbf W}^*\begin{bmatrix} V & 0 \\ 0 & V \end{bmatrix}
\begin{bmatrix}  \sqrt{H} & \sqrt{I-H} \:\\ \sqrt{I-H} & -\sqrt{H} \end{bmatrix}{\mathbf W}. \label{Ugen2}
\end{equation}
Here $U_j$, $U_{ij}$ are arbitrary unitary operators acting on ${\mathcal M}_{jj}$ and from ${\mathcal M}_{ji}$ onto ${\mathcal M}_{ij}$, respectively, and
$V$ is an arbitrary unitary operator acting on $\mathcal M$ and commuting with $H$.  \end{theorem}

{\em Proof.}
Let $\mathbf U$ be the unitary operator \eqref{U} and let $U$ be a unitary operator satisfying~\eqref{uu}.
We may write $U=Z\mathbf U$ with a unitary operator $Z$. If \eqref{uu} holds, then $Z$
commutes with both $P$ and $Q$. We may represent $Z$ in the form
\[ Z= \left(\bigoplus_{i,j=0,1}Z_{ij}\right)\oplus{\bf W}^*\begin{bmatrix} V_{11} & V_{12} \\ V_{21} & V_{22} \end{bmatrix}{\bf W}.\]
From \eqref{Halm1} and the equality $ZP=PZ$ we obtain $V_{12}=V_{21}=0$. Representation~\eqref{Halm2} and the relation $ZQ=QZ$ imply
that $V_{11}$ and $V_{22}$ commute with $H$ and that $\sqrt{H(I-H)}V_{11}=V_{22}\sqrt{H(I-H)}$. Since $V_{22}$ commutes with $H$,
it follows that $\sqrt{H(I-H)}(V_{11}-V_{22})=0$ and thus $H(I-H)(V_{11}-V_{22})=0$. As both $H$ and $I-H$ are injective,
we conclude that $V_{11}=V_{22}=:V$. Thus,
\[ Z= \left(\bigoplus_{i,j=0,1}Z_{ij}\right)\oplus{\bf W}^*\begin{bmatrix} V & 0 \\ 0 & V \end{bmatrix}{\bf W}.\]
This operator is unitary if and only if so are $Z_{ij}$ and $V$. This shows that $U=Z\bU$ has the form~\eqref{Ugen2}.
Conversely, using~\eqref{Halm1} and~\eqref{Halm2} it can be readily verified that an operator $U$ as in the theorem
satisfies~\eqref{uu}.
$\;\:\square$

\section{Additional remarks} \label{S Add}

Recall (\cite{GiKu}, see also \cite{Spit94} and \cite[Theorem 7.1]{BSpit10}) that the von Neumann algebra ${\mathcal A}(P,Q)$
generated by the pair $P,Q$ consists of all the operators of the form
\begin{equation}
\left(\bigoplus_{i,j=0,1}a_{ij}I_{{\mathcal M}_{ij}}\right)\oplus{\bf W}^*\begin{bmatrix}\phi_{00}(H) &  \phi_{01}(H) \\
\phi_{10}(H) &  \phi_{11}(H)\end{bmatrix}{\bf W}, \label{nal}
\end{equation}
where $a_{ij}\in\C$ and $\phi_{ij}$ functions on $[0,1]$ that are Borel-measurable and essentially bounded with respect to the spectral measure of $H$.
From \eqref{Ugen2} and \eqref{nal} we immediately obtain the following.

\begin{cor}\label{co:alg}
Let $P,Q$ be orthogonal projections onto $\mathcal L$ and $\mathcal N$, respectively. There exist unitary $U\in{\mathcal A}(P,Q)$
satisfying \eqref{uu} if and only if
\begin{equation}
\cM_{01}={\mathcal L}\cap{\mathcal N}^\perp=\{0\}, \quad
\cM_{10}={\mathcal N}\cap{\mathcal L}^\perp=\{0\}.\label{LN}
\end{equation}
If these two equalities hold, then all such $U$ are given by the formula
\begin{equation}
U=a_0I_{{\mathcal M}_{00}}\oplus a_1I_{{\mathcal M}_{11}}\oplus {\mathbf W}^*\begin{bmatrix} \phi(H) & 0 \\ 0 & \phi(H) \end{bmatrix}
\begin{bmatrix}  \sqrt{H} & \sqrt{I-H}\: \\ \sqrt{I-H} & -\sqrt{H} \end{bmatrix}{\mathbf W}, \label{Ualg}
\end{equation}
where $\abs{a_0}=\abs{a_1}=1$ and $\phi$ is a Borel measurable unimodular function on $[0,1]$.
\end{cor}

We pass now to the $C^*$-algebra ${\mathcal B}(P,Q)$ generated by $P$ and $Q$. Its elements are characterized among all those from ${\mathcal A}(P,Q)$ by functions $\phi_{ij}$ in \eqref{nal} that are continuous, not just measurable, on the spectrum $\Delta$ of $H$ and satisfy the following additional conditions:

\medskip
if $0 \in \Delta$ then $\phi_{01}(0)=\phi_{10}(0)=0$,

if $0 \in \Delta$ and $\cM_{00} \neq \{0\}$ then $\phi_{11}(0) =a_{00}$,

if $0 \in \Delta$ and $\cM_{11} \neq \{0\}$ then $\phi_{00}(0) = a_{11}$,

if $1 \in \Delta$ then $\phi_{01}(1) = \phi_{10}(1) = 0$,

if $1 \in \Delta$ and $\cM_{01} \neq \{0\}$ then $\phi_{11}(1) = a_{01}$,

if $1 \in \Delta$ and $\cM_{10} \neq \{0\}$ then $\phi_{00}(1) = a_{10}$;

\medskip
\noindent
see \cite[Section 4]{BSpit10} and references therein keeping in mind a slight notational deviation caused by the fact that the operator $H$ in \cite{BSpit10} is our $I-H$.
\begin{cor}\label{co:calg}Let $P,Q$ be orthogonal projections onto $\mathcal L$ and $\mathcal N$, respectively. There exist unitary $U\in{\mathcal B}(P,Q)$
satisfying \eqref{uu} if and only if the operator $P+Q-I$ is invertible.
If this is the case, then all such $U$ are given by formula \eqref{Ualg}
in which $\abs{a_0}=\abs{a_1}=1$ and $\phi$ is a continuous unimodular function on $[0,1]$.
\end{cor}

{\em Proof.}
By \eqref{Halm1} and \eqref{Halm2}, the invertibility of $P+Q-I$ is equivalent to condition \eqref{LN} combined with the invertibility of $H$. The former condition is necessary due to Corollary~\ref{co:alg}. To prove the necessity of the latter, observe that in the notation of \eqref{nal}
for the operators \eqref{Ualg} we have $\phi_{01}(t)=\phi_{10}(t)=\phi(t)\sqrt{1-t}$. Since $\phi$ is unimodular, $\phi(0)$ cannot be zero.
It follows that $\phi_{01}(0) \neq 0$ and $\phi_{10}(0) \neq 0$. Thus, by the first of the six additional conditions
listed above, $0$ is not in $\Delta$.

Conversely, let \eqref{LN} be valid and $0\notin\Delta$. Suppose $U$ is given by~\eqref{Ualg} with
$\abs{a_0}=\abs{a_1}=1$ and a continuous unimodular function $\phi$ on $[0,1]$. Passing to notation~\eqref{nal},
we have again $\phi_{01}(t)=\phi_{10}(t)=\phi(t)\sqrt{1-t}$, implying
that $\phi_{01}(1)=\phi_{10}(1)=0$. Thus, we need not take care of the first three of the six additional
requirements listed above, and the last three of them are automatically satisfied. It results that all operators of the form \eqref{Ualg} with continuous (and not
just measurable) unimodular $\phi$ belong to ${\mathcal B}(P,Q)$.
$\;\:\square$

\medskip
Finally, recall that if $H$ is a selfadjoint operator with a simple spectrum, then the only operators commuting with $H$ are functions of $H$;
see, e.g., Lemma~5.4.9 of~\cite{SimonOT}.
Combining Theorem~\ref{th:genform} for generically positioned projections with Corollary~\ref{co:alg} we therefore obtain
the following.

\begin{cor}\label{co:allin}
Let $P,Q$ be a pair of generically positioned orthogonal projections such that the spectrum of $PQP$ is simple. Then all
unitary operators $U$ satisfying \eqref{uu} lie in ${\mathcal A}(P,Q)$.
\end{cor}

Indeed, by Theorem~2 all such $U$ have the form \eqref{Ugen2}. The unitary operator $V$ there, commuting with $H$, has to be a function of the latter. Thus, $U$
is actually of the form of the last direct summand of \eqref{Ualg}. The other two summands are not present, since $P,Q$ are generically positioned.
So, $U\in{\mathcal A}(P,Q)$ due to Corollary~\ref{co:alg}.

\end{document}